\documentclass[a4paper,12pt]{article}

\newcommand{\keys}{0}       
\newcommand{\headings}{0}   
\newcommand{\comments}{0}   
\newcommand{\changes}{0}    
\newcommand{\details}{0}    

\title{Gradients on sets}

\author{ Jan Mankau, Friedemann Schuricht \\
         {\small TU Dresden, Fakult\"at Mathematik} \\ 
         {\small 01062 Dresden, Germany} 
       }

\usepackage[a4paper,left=30mm,top=32mm,right=30mm,bottom=33mm]{geometry}
\date{}
\usepackage{amssymb}
\usepackage{amsmath}
\usepackage{amsthm}
\usepackage{mathrsfs}  

\usepackage{marginnote}                
           \reversemarginpar           

\usepackage{ifthen}
\usepackage{eurosym}

\usepackage{varioref}        

\usepackage{pdfpages}

\providecommand{\includecolor}{1}            
\ifthenelse{ \equal{\includecolor}{1} }      %
  {\usepackage{color}}
  {}

\providecommand{\keys}{0}               
\ifthenelse{ \equal{\keys}{1} }         
  {\usepackage{showkeys}}{}
\ifthenelse{ \equal{\keys}{1} }
  {\newcommand{\alb}[1]{
    \ifthenelse{ \equal{#1}{} }
    {}
    {\marginnote{{\scriptsize \fbox{#1}}}\label{#1}}  }}
  {\newcommand{\alb}[1]{\label{#1}}}
\ifthenelse{ \equal{\keys}{1} }
  {\newcommand{\llb}[1]{
    \ifthenelse{ \equal{#1}{} }
    {}
    {\hspace*{-5mm}\makebox[0mm][r]{\scriptsize \fbox{#1}} \hspace*{5mm}} }}
  {\newcommand{\llb}[1]{}}
\ifthenelse{ \equal{\keys}{1} }
  {\newcommand{\lbnn}[1]{
    \ifthenelse{ \equal{#1}{} }
    {}
    {\hspace*{-5mm}\makebox[0mm][r]{\scriptsize \fbox{#1}} \hspace*{5mm}}\nn}}
  {\newcommand{\lbnn}[1]{\nn}}
\ifthenelse{ \equal{\keys}{1} }
  {\newcommand{\lbn}[1]{
    \ifthenelse{ \equal{#1}{} }
    {}
    {\hspace*{-5mm}\makebox[0mm][r]{\scriptsize \fbox{#1}} \hspace*{5mm}}
                                                              \label{#1}} }
  {\newcommand{\lbn}[1]{\label{#1}}}

\providecommand{\headings}{0}           
\ifthenelse{ \equal{\headings}{1} }     
  {\pagestyle{myheadings}
   \markboth{\footnotesize\sf File:\quad\jobname\qquad Version:\quad\today}
            {\footnotesize\sf File:\quad\jobname\qquad Version:\quad\today}}
  {}        

\ifthenelse{ \equal{\headings}{1} }     
  {\newcommand{\headingfrontpage} 
     {\renewcommand{\thefootnote}{}
      \footnotetext{{\sf File:\quad\jobname\qquad Version:\quad\today}} 
      \renewcommand{\thefootnote}{\arabic{footnote}}}}
  {\newcommand{\headingfrontpage}{}}

\providecommand{\comments}{0}        
\ifthenelse{ \equal{\comments}{1} }  
  {\newcommand{\comment}[5]
    {\textcolor{#4}{#1{#3 #5} #2}}}
  {\newcommand{\comment}[5]{}}

\providecommand{\changes}{0}            
\ifthenelse{ \equal{\changes}{1} }      
  {\newcommand{\bchange}[2] {\color{#2} #1 }}  
  {\newcommand{\bchange}[2]{}}
\ifthenelse{ \equal{\changes}{1} }      
  {\newcommand{\echange}[1]{#1\normalcolor}}
  {\newcommand{\echange}[1]{}}

\providecommand{\details}{0}            
\ifthenelse{ \equal{\details}{1} }      
  {\newcommand{\detail}[1]{\par{\footnotesize\color{blue} {\sf Detail}:
                                 #1\par\noindent}}}
  {\newcommand{\detail}[1]{}}

\numberwithin{equation}{section}      

\newtheorem{theo}[equation]{Theorem} 
\newtheorem{prop}[equation]{Proposition}

\newtheorem{cor}[equation]{Corollary}

\theoremstyle{definition} 
\theoremstyle{definition} \newtheorem{rem}[equation]{Remark}  
\theoremstyle{definition} \newtheorem{exam}[equation]{Example}  

\newcounter{nr}  \labelformat{nr}{{\normalfont \arabic{nr}}}

\newcommand{\bgl}[1][{\normalfont (\arabic{nr})}]
   {\begin{list} {#1} {\usecounter{nr}
          \setlength{\topsep}{0.5ex plus0.2ex minus0.1ex}     
          \setlength{\itemsep}{0.2ex plus0.05ex minus0.03ex}  
          \parsep0pt \itemindent0pt 
          \leftmargin30pt   \labelwidth20pt}   }
\newcommand{\el}{\end{list}}

\newcounter{nra}  \labelformat{nra}{{\normalfont \alph{nr}}}
\newcommand{\bgla}[1][{\normalfont (\alph{nra})}]
   {\begin{list} {#1} {\usecounter{nra}
          \setlength{\topsep}{0.5ex plus0.2ex minus0.1ex}     
          \setlength{\itemsep}{0.2ex plus0.05ex minus0.03ex}  
          \parsep0pt \itemindent0pt 
          \leftmargin30pt   \labelwidth20pt}   }
\newcommand{\ela}{\end{list}}

\newcommand{\bgln}                
  {\bg{list}{(\arabic{section}.\arabic{equation})}{\usecounter{equation} 
          \setlength{\topsep}{0.5ex plus0.2ex minus0.1ex}     
          \setlength{\itemsep}{0.2ex plus0.05ex minus0.03ex}  
          \parsep0pt \itemindent0pt 
          \leftmargin30pt   \labelwidth20pt}   }

\newcommand{\beq}{\begin{equation}\nonumber}
\newcommand{\bee}{\begin{equation}\nonumber}
\newcommand{\bn}[1]{\begin{equation}\label{#1}} 
\newcommand{\bnn}[1]{\begin{equation}\llb{#1}\nonumber}

\newcommand{\ee}{\end{equation}}

\newcommand{\ba}{\begin{eqnarray}} \newcommand{\ea}{\end{eqnarray}}
\newcommand{\barr}{\begin{array}} \newcommand{\earr}{\end{array}}

\newcommand{\lb}{\label}
\newcommand{\nn}{\nonumber}

\newcommand{\ben}[1][]{\begin{equation}\label{#1}}
\newcommand{\benn}[1][]{\begin{equation}\llb{#1}\nonumber}

\newcommand{\en}{\end{equation}}

\newenvironment{pf}[1]{\noindent {\sc Proof}#1.\hspace{1pt}}%
     {\hspace*{\fill} $\diamondsuit$\par\medskip}

\newcommand{\bg}[1]{\begin{#1}} \newcommand{\e}[1]{\end{#1}}

\newcommand{\z}{\enspace}

\newcommand{\qmz}[1]{\quad\mbox{#1}\z}
\newcommand{\qmq}[1]{\quad\mbox{#1}\quad}
\newcommand{\zmz}[1]{\z\mbox{#1}\z}
\newcommand{\leer}{\mbox{}}
\newcommand{\tx}[1]{\text{#1}}

\newcommand{\msk}{\medskip}
   \newcommand{\bbsk}{\bigskip\bigskip}

\newcommand{\noi}{\noindent}
  \newcommand{\ol}{\overline}

\newcommand{\reff}[1]{{\rm (\ref{#1})}}
\newcommand{\refr}[1]{{\rm \ref{#1}}}

\newcommand{\bmp}[1]{\begin{minipage}[t]{#1}}
\newcommand{\bmpc}[1]{\begin{minipage}[c]{#1}}
\newcommand{\emp}{\end{minipage}}

\newcommand{\op}[1]{\operatorname{#1}}   

\newcommand{\osm}[2]{\ensuremath{\overset{#1}{#2}}}

\newcommand{\wto}{\rightharpoonup}
\newcommand{\wsto}{\mathop{\wto}\limits^*}

\newcommand{\dl}{\langle}  \newcommand{\dr}{\rangle}

\newcommand{\pa}{\partial}

\newcommand{\subs}{\subset}

\renewcommand{\d}[1][]{{\rm d}}     

\newcommand{\ti}{\tilde}

\newcommand{\N}{\ensuremath{\mathbb{N}}} 
 
\newcommand{\R}{\ensuremath{\mathbb{R}}}

\def\XXint#1#2#3{{\setbox0=\hbox{$#1{#2#3}{\int}$}
\vcenter{\hbox{$#2#3$}}\kern-.5\wd0}}

\newcommand{\al}{\alpha}

\newcommand{\de}{\delta}      
\newcommand{\ep}{\varepsilon}

\newcommand{\la}{\lambda}     

\newcommand{\si}{\sigma}      
\newcommand{\ta}{\tau}

\newcommand{\gG}[3]{\partial #1(#2)}
\newcommand{\conv}{ \,} 
\newcommand{\gdD}[4]{#1 ^0_{#2}( #3; #4) }

\newcommand{\DefFunktion}[3]{#1:\,#2\to #3 }

\newcommand{\tme}{\:\big|\:}    
\newcommand{\sk}[3]{\langle #1, #2\rangle} 
\newcommand{\norm}[1]{\| #1 \| }       
\newcommand{\norms}[2]{\left\| #1 \right\|_{#2} }       
\newcommand{\set}[2]{\left\{ #1 \tme #2\right\} }               
\newcommand{\sset}[1]{\left\{ #1 \right\} }               
\newcommand{\ball}[3]{ B_{#2}(#1) }         

\newcommand{\cball}[3]{\overline {B_{#2}(#1)}} 

\usepackage{graphicx}

\allowdisplaybreaks

\begin{document}

\bibliographystyle{plain}


\maketitle
\headingfrontpage               


 
\vspace{10mm}

\bg{abstract} For a locally Lipschitz continuous function $f:X\to\R$  
the generalized gradient $\pa f(x)$ of Clarke is used to develop 
some (set-valued) gradient on a set $A\subs X$. Existence, uniqueness and some
approximation are considered for optimal descent directions on set $A$.
The results serve as basis for nonsmooth numerical descent algorithms that can
be found in subsequent papers. 
\e{abstract}

\vspace{4mm}

\section{Introduction}

For a smooth function $f:X\to\R$ the derivative $f'(x)$ in particular indicates
directions of descent near $x$. This fact serves as basis for typical
numerical descent algorithms. However such algorithms fail in cases where the
direction of descent changes rapidly in a small neighborhood of $x$. 
This typically occurs
for functions having large second derivatives and, even worse, for functions   
that are not differentiable. In such
situations it becomes necessary to use more information of $f$ for the
selection of a descent direction. If we consider some $f$ being the pointwise 
maximum of two (non-constant) linear functions, we have to realize that also 
Clarke's set-valued generalized gradient $\pa f(x)$, defined for Lipschitz 
continuous functions $f$,
does not provide enough information for a stable scheme. 
Therefore the selection of a robust descent direction is only possible 
if one uses relevant information of $f$ from some suitable neighborhood of
$x$. 

We consider locally Lipschitz continuous functions
$f:X\to\R$ on a Banach space~$X$. 
Using the generalized gradients of Clarke we introduce    
some (set-valued) gradient $\pa f(A)$ of $f$ on a set $A\subs X$ and, 
with Clarke's generalized directional derivative $\gdD{f}{}{y}{h}$, 
we define some directional derivative $\gdD{f}{}{A}{h}$ of $f$ on $A$ in
direction~$h$. In Section~\ref{gg} we verify basic properties 
for these new quantities where some are quite similar to that in
Clarke's calculus. For sequences
of sets $A_k\to A$ converging in the Hausdorff metric,
some general upper semicontinuity is shown. Here the relevance of certain
assumptions is illuminated by examples. Moreover we show that 
the $\ep$-generalized gradient $\de_\ep f(x)$ of $f$ at $x$
introduced in Goldstein \cite{Goldstein} for $X=\R^n$ (that somehow relies
on Rademacher's Theorem for Lipschitz continuous functions)
agrees with $\gG{f}{\cball{x}{\ep}{}}{\conv}$. Finally 
we consider regularity in the sense that $0\not\in\pa f(A)$ and,
in particular, a result from Goldstein \cite{Goldstein} is extended 
to Banach spaces. 
In Section \ref{odd} we define descent directions and optimal descent
directions of $f$ on $A$ . Then existence and general properties of optimal
descent directions are analyzed. An example demonstrates that there might be
no optimal descent direction in a non-reflexive Banach space.    
Uniqueness of an optimal descent direction can be verified 
for strictly convex Banach spaces. Examples show that the selection of 
descent directions and optimal descent directions needs much more care 
in spaces that are not strictly convex. Furthermore we provide 
some stability and approximation results for optimal descent directions 
that are very useful for applications in
numerics. The advantage of gradients on sets and corresponding descent
directions for numerical algorithms is 
demonstrated by a simple but typical example.
Applications of the analytical results to nonsmooth descent algorithms and 
corresponding numerical simulations can be found in subsequent papers. 

\msk

{\it Notation:} By $X$ we denote a Banach space, by $X^*$ its dual, and 
by $\dl\cdot,\cdot\dr$ the corresponding duality pairing.
We call $X$ (or $X^*$) strictly or uniformly convex if the norm has that 
property (cf. \cite{Cioranescu}).
For a set $M$ we use $\ol{M}$ for its closure,
$\op{conv}M$ for its convex hull, and  
$\overline{\op{conv}}^*M$ for its weak$^*$-closed convex hull.
$B_\ep(x)$ stands for the open $\ep$-neighborhood of point $x$ and  
$B_\ep(M)$ for the open $\ep$-neighborhood of set $M$.
We write $]x,y[$ and $[x,y]$ for the open and closed segment (or interval),
respectively, generated by the points $x,y$. 
Clarke's generalized directional derivative is denoted by $f^0(x;h)$ and its 
generalized gradient by $\pa f(x)\subs X^*$ (cf. Clarke~\cite{Clarke}).
Notice that $\pa f(A)$ denotes the gradient defined in \reff{gg-e1} and does
not mean $\bigcup_{x\in A}\pa f(x)$.

\section{Gradients on sets}
\label{gg}

Let $X$ be a Banach space 
and let $\DefFunktion{f}{X}{\R}$
be a locally Lipschitz continuous function.
We denote the generalized gradient at $x$ by $\gG{f}{x}{}$
and the generalized directional derivative at $x$ in direction $y$
by $\gdD{f}{}{x}{y}$ (cf. Clarke~\cite{Clarke}). 
While these quantities somehow express the behavior of $f$ at the point $x$,
we are interested in information that expresses the behavior of $f$ on a whole
set. Therefore we introduce some set-valued gradient of $f$ on a set 
$A\subs X$ by using Clarke's pointwise quantities. 
Later sets $A=\cball{x}{\ep}{}$ with $\ep>0$ will be of particular interest. 

For $A\subset X$ ($\ne\emptyset$) we define the {\it gradient} of $f$ on $A$
by 
\begin{equation}\lb{gg-e1}
 \gG{f}{A}{\conv}:=\overline{\op{conv}}^*\bigcup\limits_{y\in A } \gG{f}{y}{} 
\end{equation}
(where $\overline{\op{conv}}^*$ denotes the weak$^*$ closure of the convex
hull)  
and the {\it directional derivative} of $f$ at $A$ in direction $h\in X$ by
\begin{equation}
 \gdD{f}{}{A}{h}:=\sup\limits_{y\in A }\gdD{f}{}{y}{h}\,.
\end{equation}
Clearly $\pa f(x)=\pa f(\{x\})$ and 
$\gdD{f}{}{x}{h}=\gdD{f}{}{\{x\}}{h}$. 
Let us start with some basic properties. 
\begin{prop}\label{gg-s1}
Let $A\subs X$ be nonempty and let 
$\DefFunktion{f}{X}{\R}$ be Lipschitz continuous of rank~$L$ on a
neighborhood of $A$. Then: 
\bgl
 \item $\gG{f}{A}{\conv}$ is nonempty, convex, weak$^*$-compact and bounded by
   $L$. 
 \item $\gdD{f}{}{A}{\cdot}$ is finite, positively homogeneous, subadditive,
   and Lipschitz continuous of rank $L$. Moreover it is the support function
   of $\gG{f}{A}{\conv}$ with 
\ben
\gdD{f}{}{A}{h} = \max\limits_{a\in\gG{f}{A}{\conv}} 
 \sk{a}{h}{}  \qmz{for all} h\in X\,. \label{supportgG}
\ee
 \item We have 
\ben
  \gG{f}{A}{\conv} = \set{a\in X^*}{\sk{a}{h}{}\leq  \gdD{f}{}{A}{h}\zmz{for
      all} h\in X }\,.\label{ChargG}
\ee
\item \label{gg-s1-4}
  Let $h\in X$ with $\gdD{f}{}{A}{h}<0$, let $x\in A$, and let $t>0$
  with $]x,x+th[\:\subset A$. Then 
$$f(x+th)\leq f(x)+t\gdD{f}{}{A}{h}<f(x)\; .$$
\el 
\end{prop}
\begin{pf}{}
For (1) we recall that $\gG{f}{y}{}$ is nonempty and
    bounded by $L$ for all $y\in A$ (cf. \cite[Prop.~$2.1.2$]{Clarke}). 
    Thus the stated properties follow easily from the  
    definition of $\pa f(A)$ and the Banach Alaoglu Theorem. 

For (2) we first notice that $\gdD{f}{}{y}{\cdot}$ is the support
function of $\gG{f}{y}{}$ (cf. \cite[Prop.~$2.1.2$]{Clarke}). Therefore we
    obtain for the support function of
    $\gG{f}{A}{\conv}$ at $h\in X$  
 \begin{eqnarray*}
  \sup\limits_{a\in\gG{f}{A}{\conv}}     \sk{a}{h}{} 
&=&  \sup\bigg\{ \sk{a}{h}{} \,\Big|\: 
a\in \, \ol{\op{conv}}^*\bigg( \bigcup_{y\in A } \gG{f}{y}{} \bigg) \bigg\} \\
&=&  \sup\bigg\{ \sk{a}{h}{} \,\Big|\:  a\in \,\op{conv}
\bigg( \bigcup\limits_{y\in A } \gG{f}{y}{} \bigg) \bigg\}     \\
&=&  \sup\bigg\{ \sk{a}{h}{} \,\Big|\: a\in\,\bigg( \bigcup_{y\in A} 
\gG{f}{y}{} \bigg) \bigg\}     \\
&=& \sup_{y\in A} \sup_{a\in \gG{f}{y}{} } \sk{a}{h}{} 
= \:\sup_{y\in A} \gdD{f}{}{y}{h }
\, = \: \gdD{f}{}{A}{h}\,.
 \end{eqnarray*}
Since $\gG{f}{A}{\conv}$ is weak$^*$-compact, the supremum is attained and
(\ref{supportgG}) follows. The remaining properties are now easy consequences. 

For (3) we notice that characterization (\refr{ChargG}) is as general
property of support functions (cf. \cite[Prop.~$2.1.4$]{Clarke}).

For (4) we use Lebourg's mean value theorem 
(cf.~\cite[Prop.~$2.3.7$]{Clarke}) to get some $z\in]x,x+th[$ and 
some $a\in \gG{f}{z}{}\subset
\gG{f}{A}{\conv}$ such that 
\begin{equation*}
   f(x+th)-f(x)=\sk{a}{th}{}\stackrel{(\ref{supportgG})}{\leq} t
   \gdD{f}{}{A}{h}<0 \,,
\end{equation*}
which directly implies the assertion.
\end{pf}
\begin{prop}[upper semicontinuity]\label{prop-usc}
Let $\DefFunktion{f}{X}{\R}$ be locally Lipschitz continuous, let 
$h_k\to h$ in $X$, and let $A_k,A\subset X$ with $A$ compact and 
$A_k\to A$ in the Hausdorff metric, i.e. 
\begin{equation*}
d(A_k,A):= \inf\set{\de>0}{A\subset B_\de(A_k)  \zmz{and} 
  A_k\subset B_\de(A) }\xrightarrow{k\to\infty} 0  \,. 
\end{equation*}
Then
\begin{eqnarray}
 \limsup_{k\to\infty}\gdD{f}{}{A_k}{h_k}  
&\leq&  \gdD{f}{}{A}{h}\label{usc-dirctional}\\
\big\{ a\in X^* \mid a_k\wsto a \zmz{for} a_k\in\gG{f}{A_k}{\conv} \big\}  
&\subs &\gG{f}{A}{\conv}\,.\label{usc-grad}
\end{eqnarray}
If $A\subs A_k$ for all $k\in\N$, then we have equality in 
\reff{usc-grad} and
\begin{equation}
 \lim_{k\to\infty}\gdD{f}{}{A_k}{h_k} =  \gdD{f}{}{A}{h}\label{usc-dirctionaleq}
\end{equation}
\end{prop}
\noi
With $A=\{x\}$ we directly derive the following statement. 
\begin{cor}
Let $x\in X$ and $\ep_k\to0$ such that 
$x\in A_k\subs B_{\ep_k}(x)$ and let $h\in X$. Then 
\begin{equation*}
 \lim_{k\to\infty}\gdD{f}{}{A_k}{h}=\gdD{f}{}{x}{h}
 \qmq{and} \bigcap_{k\in\N}\gG{f}{A_k}{\conv}=\gG{f}{x}{}\,.  
\end{equation*}
\end{cor}

\begin{pf}{ of Proposition \ref{prop-usc}}
By definition and assumption there exist $x_k\in A_k$ and $z_k\in A$ with
\bee
 \gdD{f}{}{A_k}{h_k}\geq \gdD{f}{}{x_k}{h_k}\geq \gdD{f}{}{A_k}{h_k}-\frac1k \qmq{and}
 \norm{x_k-z_k}\to 0\,.
\ee
By compactness of $A$ we get, possibly for a subsequence, 
\bee
  \gdD{f}{}{x_k}{h_k}\to\limsup_{k\to\infty}\gdD{f}{}{A_k}{h_k}
  \qmq{and} z_k\to :z\in A\,.
\ee
Consequently $x_k\to z$. Since $\gdD{f}{}{\cdot}{\cdot}$ is upper
semicontinuous (cf. \cite[Prop.~$2.1.1$]{Clarke}),   
\begin{equation*}
 \gdD{f}{}{A}{h}\geq\gdD{f}{}{z}{h}\geq \lim_{k\to\infty}\gdD{f}{}{x_k}{h_k}=
 \limsup_{k\to\infty}\gdD{f}{}{A_k}{h_k}
\end{equation*}
and we have (\ref{usc-dirctional}). 

Let now $a_k\in\gG{f}{A_k}{\conv}$ with
$a_k\wsto a$. Hence  
\begin{equation*}
 \gdD{f}{}{A}{h} \geq
 \limsup_{k\to\infty}\gdD{f}{}{A_k}{h}\stackrel{}{\geq}
 \lim_{k\to\infty} \sk{a_k}{h}{}=\sk{a}{h}{} \qmz{for all} h\in X\,. 
\end{equation*}
Thus $a\in \gG{f}{A}{\conv}$ by (\ref{ChargG}). If $A\subs A_k$, then
$\gdD{f}{}{A}{\cdot}\leq \gdD{f}{}{A_k}{\cdot}$ and $\gG{f}{A}{\conv}\subs
\gG{f}{A_k}{\conv}$ by definition. Hence, equality in (\ref{usc-grad}) follows
in the case that $A\subs A_k$ for all $k\in\N$. Furthermore
\begin{equation*}
\gdD{f}{}{A}{h} \leq    \liminf\limits_{k\to\infty}\gdD{f}{}{A_k}{h_k}  \leq
\limsup\limits_{k\to\infty}\gdD{f}{}{A_k}{h_k}
\stackrel{(\ref{usc-dirctional})}{\leq } \gdD{f}{}{A}{h} 
\end{equation*}
and (\ref{usc-dirctionaleq}) follows.
\end{pf}

\begin{exam} We present some examples showing the necessity of central
  assumptions in Proposition~\ref{prop-usc}.

\bgl
 \item Let $f:\R\to\R$ be given by $f(x)=|x|$ and let 
\bee
  h_k=h=-1\,,\quad A_k:=]-\tfrac1k,1[\,,\quad A:=]0,1[\,.
\ee 
Obviously $d(A_k,A)\to 0$, but $A$ is not compact. 
We have $\gG{f}{A}{\conv}=\sset{1}$ and
\begin{equation*}
 \gdD{f}{}{A_k}{h_k}\geq \gdD{f}{}{0}{h}=1>-1\stackrel{(\ref{supportgG})}{=}
 \gdD{f}{}{A}{h}\,. 
\end{equation*}
Hence \reff{usc-dirctional} is not satisfied. 

\item Let again $f:\R\to\R$ be given by $f(x)=|x|$ and let 
\bee
  h_k=h=-1\,,\quad A_k:=[\tfrac1k,1]\,,\quad A:=[0,1]\,.
\ee
Here $A$ is compact and $d(A_k,A)\to0$, but $A\not\subs A_k$. We have 
$\gG{f}{A_k}{\conv}=\sset{1}$ and
\begin{equation*}
 \gdD{f}{}{A}{h} \geq \gdD{f}{}{0}{h}=1 > -1\stackrel{(\ref{supportgG})}{=}
 \gdD{f}{}{A_k}{h_k}\,. 
\end{equation*}
Therefore (\ref{usc-dirctional}) is satisfied, but without  
equality as in \reff{usc-dirctionaleq}. 

\item For $X=\ell^2$ (sequences $x=(\xi_i)_{i\in\N}$ in $\R$ with
$\norm{x}^2=\sum_{k\in\N}|\xi_i|^2<\infty$) we consider
\bee
  A:=\sset{0}\,,\quad
  A_k:=\cball{0}{1}{}\cap \set{(\xi_i)\in\ell^2}{\xi_j=0\;\mbox{for}\;j<k}.
\ee
Clearly $A$ is compact and $A=\bigcap\limits_{k\in\N}A_k$.
But $d(A_k,A)\not\to0$, since there are $x_k\in A_k$ with $\norm{x_k}=1$.
With fixed $z\in\ell^2\setminus\{0\}$ and $\phi\in C^{\infty}(\R,\R)$ 
satisfying
\bee
  \phi(\al)=1 \zmz{for} \al<\tfrac14\,, \quad
  \phi(\al)=0 \zmz{for} \al>\tfrac34\,, 
\ee
we define $f:\ell^2\to\R$ by  
\bee
  f(x):=\sk{z}{x}{}\,\phi(\norm{x})\,.
\ee
Obviously $f$ is locally Lipschitz continuous with
\bee
  \gdD{f}{}{x}{-z}=0 \zmz{if} \norm{x}=1 \qmq{and}
  \gdD{f}{}{0}{-z}=-\norm{z}^2\neq 0\,.
\ee 
For $h_k=h=-z$ we obtain
\begin{equation*}
 \gdD{f}{}{A_k}{h_k}\geq \gdD{f}{}{x_k}{-z}=0>-\norm{z}^2=
 \gdD{f}{}{0}{-z}= \gdD{f}{}{A}{h}
\end{equation*}
and, again, (\ref{usc-dirctional}) is violated.
\el
\end{exam}

\bbsk

For $X=\R^n$ and $\ep\ge 0$ the {\it $\ep$-generalized gradient} 
of $f$ at $x\in X$ is given according to 
Goldstein \cite{Goldstein} by
\begin{equation}\label{gg-e2}
 \de_{\ep}f(x):=\op{conv} \bigcap\limits_{k=1}^{\infty}
  \overline{\set{f'(y)}{
y\in\cball{x}{\ep+\frac{1}{k}}{},\; f'(y)\;\mbox{exists}}} \,
\end{equation}
(where $f'(x)$ denotes the usual derivative).
\begin{cor}\label{gg-s4}
Let $f:\R^n\to\R$ be locally Lipschitz continuous. Then
\begin{equation*}
 \de_{\ep}f(x)=\gG{f}{\cball{x}{\ep}{}}{\conv} 
 \qmz{for all} x\in \R^n,\z \ep\geq0\,.
\end{equation*}
\end{cor}
\begin{pf}{}
Using the characterization of $\gG{f}{x}{}$ in $\R^n$ 
(cf. \cite[Theorem~$2.5.1$]{Clarke}), we get
\begin{equation*}
 \gG{f}{\cball{x}{\ep}{}}{\conv} \subs
 \de_{\ep}f(x) \osm{\reff{gg-e2}}{\subs} 
 \bigcap_{k\in\N} \gG{f}{\cball{x}{\ep+\frac1k}{}}{\conv}
 \osm{\reff{usc-grad}}{\subs}
 \gG{f}{\cball{x}{\ep}{}}{\conv}\,
\end{equation*}
(most right inclusion is already an equality by 
$\cball{x}{\ep}{}\subs \cball{x}{\ep+\frac1k}{}$ for all $k\in\N$).
\end{pf}
The following statement somehow generalizes Goldstein 
\cite[Propostion~$2.8$]{Goldstein} from $X=\R^n$ to a general Banach space $X$. 
\begin{prop}\label{gg-s5}
Let $f:X\to\R$ be locally Lipschitz continuous and let
$A\subset X$ be compact such that $0\notin\gG{f}{x}{}$ 
for all $x\in A$. 
Then there exists $\ep>0$ and $\si>0$ such that 
\begin{equation}\label{gg-s5-1}
 \min\big\{ \norm{a}\;\big|\; a\in\gG{f}{\cball{x}{\ep}{}}{\conv} \big\} \geq
 \si \qmz{for all} x\in A\,.
\end{equation}
\end{prop}

\begin{pf}{}
Notice that there is a minimum in \reff{gg-s5-1}, since the 
norm $\norm{\cdot}$ in $X^*$ is weak$^*$ lower semicontinuous and
$\gG{f}{\cball{x}{\ep}{}}{\conv}$ is weak$^*$ compact. 
If the statement would be wrong, then there are
$x_k\in A$ with
\bee
  \min\big\{ \norm{a} \;\big|\; a\in\gG{f}{\cball{x_k}{\frac1k}{}}{\conv} 
  \big\} < \tfrac{1}{k} \qmz{for all} k\in\N\,. 
\ee
By compactness of $A$ we can assume that $x_k\to:x\in A$. 
Moreover we find $a_k\in\gG{f}{\cball{x_k}{\frac1k}{}}{\conv}$ with $a_k\to 0$. 
Since $d\big(\cball{x_k}{\frac1k}{},\{x\}\big)\to0$,  
Proposition~\ref{prop-usc} gives the contradiction
$0\in\gG{f}{x}{}$.
\end{pf}
Let us finally show that $0\notin\gG{f}{x}{}$ implies some regularity also in
a small neighborhood of~$x$. 
\begin{prop}\label{boundedawayzero}
Let $f:X\to\R$ be locally Lipschitz continuous and let  
$0\notin\gG{f}{x}{}$ for some $x\in X$. 
Then there exist $\ep>0$ and $h\in X$ with $\norm{h}{}=1$
such that 
\begin{equation}\label{gg-s6-1}
  -\norms{a}{}\leq \sk{ a}{ h}{}\stackrel{(\ref{supportgG})}{\leq}
  \gdD{f}{}{A}{h}<0 \qmz{for all}  A\subset \ball{x}{\ep}{}, \; 
  a\in\gG{f}{A}{\conv}      \,. 
\end{equation}
\end{prop}
\begin{pf}{} 
By $0\notin\gG{f}{x}{}$, property \reff{ChargG} with $A=\sset{x}$
provides the existence of some $h\in X$ with $\norm{h}{}=1$ and
$\gdD{f}{}{x}{h}<0$. Proposition~\ref{prop-usc} implies 
\bee
  \lim_{k\to\infty}\gdD{f}{}{\ball{x}{\frac1k}{}}{h}=
  \gdD{f}{}{x}{h}<0\,.
\ee
Hence we get the most right inequality in \reff{gg-s6-1} 
for some $\ep>0$ sufficiently small. With \reff{supportgG} we obtain for any 
$a\in\gG{f}{A}{\conv}$
\bee
  -\norms{a}{}\le \sk{ a}{ h}{}\leq\gdD{f}{}{A}{h}
\ee
which verifies the assertion.
\end{pf}

\section{Optimal descent directions}
\label{odd}

Motivated by Proposition \ref{gg-s1} \reff{gg-s1-4}
we say that $h\in X$ is a {\it descent direction} of $f$ on $A$
if $\gdD{f}{}{A}{h}<0$
(cf. also Clarke \cite[Ex. 10.7]{ClarkeA}).
We call $\ti h\in X$ {\it steepest} or {\it optimal} descent
direction of $f$ on $A$ with respect to $\norm{\cdot}$ if
\begin{equation}
  \|\ti h\|=1 \qmq{and}
  \gdD{f}{}{A}{\ti h}=\min\limits_{\norm{h}{}\leq1}  \gdD{f}{}{A}{h}<0\,.
\end{equation}
For reflexive Banach spaces the existence of optimal descent directions
follows from duality theory.  
\begin{prop}[existence of optimal descent
  directions]\label{odd-s1} 
Let $A\subs X$ be nonempty and let $f:X\to\R$ be Lipschitz continuous on a
neighborhood of $A$. Then:
\bgl
 \item There is some $\ti a\in \gG{f}{A}{\conv}$ such that
\begin{equation}\label{MinMaxvertauschung}
        \inf\limits_{\norm{h}\leq
          1}\gdD{f}{}{A}{h}=-\min\limits_{a\in\gG{f}{A}{\conv}}\norm{a} =
        -\|\ti a\|\,. 
\end{equation}
\item For every pair $(\ti a,\ti h)\in\gG{f}{A}{\conv}\times\cball{0}{1}{X}$
with 
\begin{equation}\label{NodwendigProp}
 \norm{\ti a}=\min_{a\in\gG{f}{A}{\conv}}\norm{a} \qmq{and}
 \gdD{f}{}{A}{\ti h}=\min_{\norm{h}\leq1}\gdD{f}{}{A}{h}
\end{equation}
we have
\begin{equation}\label{Glamxm}
 -\norm{\ti a}=\sk{\ti a}{\ti h}{}=\gdD{f}{}{A}{\ti h}\,.
\end{equation}
\item If $X$ is reflexive, then there exists a pair $(\ti a,\ti h)
\in\gG{f}{A}{\conv}\times\cball{0}{1}{X}$ satisfying \reff{NodwendigProp}.
\el
\end{prop}
Before providing the proof we still formulate a simple consequence. 
\begin{cor}
Let $A\subs X$ be nonempty and
let $f:X\to\R$ be Lipschitz continuous on a neighborhood of $A$. Then
\begin{equation}\label{Char0ingG}
  \inf\limits_{\norm{h}\leq1}\gdD{f}{}{A}{h}<0 \qquad\Longleftrightarrow\qquad
  0\notin\gG{f}{A}{\conv} \,.
\end{equation}
Moreover, if $0\notin\gG{f}{A}{\conv}$ and 
$(\ti a,\ti h)\in\gG{f}{A}{\conv}\times\cball{0}{1}{X}$ 
satisfies \reff{NodwendigProp}, then $\ti h$
is an optimal descent direction of $f$ on $A$.
\end{cor}
\begin{pf}{ of Propostion~\ref{odd-s1}} 
For (1) we readily see that we have a minimizer $\ti a\in\gG{f}{A}{\conv}$
and we use \reff{supportgG} to get
\begin{equation*}
 \inf_{\norm{h}{}\leq 1} \gdD{f}{}{A}{h} = 
 \inf_{\norm{h}{}\leq 1}\max_{a\in\gG{f}{A}{\conv}}\sk{ a}{h}{}\,. 
\end{equation*}
Since $\gG{f}{A}{\conv}$ is weak* compact, we can exchange inf and max
by Aubin's lopsided minimax theorem (cf. 
\cite[Theorem~$6.2.7$]{AubinEkeland1984}) and obtain
\begin{equation}\label{MinMaxProblem}
 \inf_{\norm{h}{}\leq 1} \gdD{f}{}{A}{h} = 
 \max_{a\in\gG{f}{A}{\conv} }\inf_{\norm{h}{}\leq 1}\sk{ a}{h}{} =
 \max_{a\in\gG{f}{A}{\conv}}-\norm{a} = 
 -\min_{a\in\gG{f}{A}{\conv}}\norm{a} \,.
\end{equation}

For (2) let  
$(\ti a,\ti h)\in\gG{f}{A}{\conv}\times\cball{0}{1}{X}$ satisfy
(\ref{NodwendigProp}). Then
\begin{equation*}
  -\norm{\ti a} \osm{\reff{MinMaxvertauschung}}{=}
 \gdD{f}{}{A}{\ti
   h}\stackrel{(\ref{supportgG})}{=}\max\limits_{a\in\gG{f}{A}{\conv}}\sk{a}{\ti
   h}{}\geq\sk{\ti a}{\ti h}{}\geq\inf\limits_{\norm{h}\leq 1}\sk{\ti
   a}{h}{}=-\norm{\ti a}
\end{equation*}
which readily gives \reff{Glamxm}.

For (3) we first observe that there is a minimizer 
$\ti a\in\gG{f}{A}{\conv}$ satisfying the left part of 
\reff{NodwendigProp} (cf. also (1)). For the right part we use that 
$\gdD{f}{}{A}{\cdot}$ is convex and continuous and, thus, weakly lower
semicontinuous. Since $X$ is reflexive, there is a minimizer $\ti h$ on the
bounded set $\ol{B_1(0)}$ by the Weierstra\ss\ Theorem.
\end{pf}
The following example shows that there might not be an optimal descent
direction in a non-reflexive Banach space $X$. 

\begin{exam}
For $X=c_0$ (sequences $x=(\xi_i)_{i\in\N}$ in $\R$ with $\xi_i\to 0$ and 
$\norm{x}=\max_{i\in\N}|\xi_i|$) the dual is $X^*=\ell^1$
(sequences $x=(\xi_i)_{i\in\N}$ in $\R$ with
$\norm{x}=\sum_{k\in\N}|\xi_i|<\infty$,
cf. \cite[Satz II.2.3]{Werner2005}).
Then $f=\big(\frac{1}{2^{i+1}}\big)_{i\in\N}\in X^*$
is a Lipschitz continuous function on $c_0$ with
\begin{equation*}
f(x)=\dl f,x\dr=
\sum_{i\in\N}\frac {\xi_i}{2^{i+1}}\qmq{and}\norm{ f}=1\,.
\end{equation*}
By linearity, $\gG{f}{x}{}=\gG{f}{A}{\conv}=\sset{f}$ for all $x\in c_0$ and
all nonempty $A\subs c_0$. Hence $\ti a=f$ always satisfies 
\reff{NodwendigProp} and we have
\bee
  \inf_{\|h\|\le 1}\gdD{f}{}{A}{h}=\inf_{\|h\|\le 1}\dl f,h \dr = -1 
  \qmz{for all nonempty} A\subs c_0\,.
\ee
But there is no $\ti h\in c_0$ with $\norm{\ti h}\leq 1$ such that 
$\gdD{f}{}{A}{\ti h}=-1$, i.e. there is no optimal descent direction.  
We merely find arbitrarily good approximations as
e.g. $h_k=(\xi_i^k)_{i\in\N}\in c_0$ with 
\bee
  \xi_i^k=-1 \zmz{for} i\le k\,, \quad \xi_i^k=0 \zmz{for} i>k\,.
\ee 
Obviously $\norm{h_k}=1$ and, using \reff{supportgG}, we readily get
$\gdD{f}{}{A}{h_k}=\sk{f}{h_k}{}\to-1$.
\end{exam}
\begin{theo}[uniqueness of optimal descent direction]\label{Satz1Umg}
Let $X$ be reflexive and let $X$, $X^*$ be strictly convex, let 
 $A\subs X$ be nonempty, and let $f: X\to\R$ be Lipschitz
continuous on a neighborhood of $A$. Then there is a
unique $\ti a\in\gG{f}{A}{\conv}$ with
\begin{equation}\label{Normkleinstes}
 \|\ti a\|=\min_{a\in\gG{f}{A}{\conv}} \norm{a} \,.
\end{equation}
Moreover, if $0\notin\gG{f}{A}{\conv}$, 
then there exists a unique optimal descent
direction $\ti h$ of $f$ which is characterized by
\begin{equation}\label{NormkleinstesDual}
 \dl\ti a,\ti h\dr = -\|\ti a\| \qmq{with} \|\ti h\|=1\,.
\end{equation}
\end{theo}
\begin{pf}{}
By Propostion~\ref{odd-s1} there are $\ti a$ and $\ti h$ satisfying
\reff{Normkleinstes}, \reff{NormkleinstesDual}. 
Since $X^*$ is strictly convex and $\gG{f}{A}{\conv}$ convex, 
$\ti a$ in \reff{Normkleinstes} is unique. Since $X$ is strictly convex, 
$\ti h$ in \reff{NormkleinstesDual} is also unique.
\end{pf}

\begin{rem}\leer
Notice that for every reflexive Banach space $X$ there exists an
equivalent norm such that $X$ and $X^*$ are strictly convex
(cf. \cite[Theorem~III.2.9]{Cioranescu}). However, since the 
optimal descent direction $\ti h$ depends on the norm in general, 
$\ti h$ might change by a change of norm. In particular, 
the derivative $f'(x)$ of a smooth function $f$ is independent of an
equivalent norm, but the optimal descent direction $\ti h$ on $A=\{x\}$ 
might be different for an equivalent norm.   
\end{rem}
%
%
%

\begin{exam}
We consider $X:=\R^2$ with the non strictly convex norms
$\norm{x}_1$ (1-norm) and $\norm{x}_\infty$ (maximum norm). 
We will demonstrate that the selection of a descent direction needs more care
in a reflexive but not strictly convex space where 
\reff{NormkleinstesDual} is not sufficient for the selection.
\bgl
\item Let $X=(\R^2,\norm{\cdot}_1)$ and, thus, its dual
$X^*=(\R^2,\norm{\cdot}_\infty)$. We define $f:\R^2\to\R$ by 
\bee
  f(x_1,x_2)=x_1+|x_2| \,.
\ee
With $x=(1,0)$ and $A=\sset{x}$ we get
\bee
  \gG{f}{A}{\conv}=\gG{f}{x}{}=\{(1,\la) \mid \la\in[-1,1]\}\,.
\ee
Obviously any $\ti a \in\gG{f}{A}{\conv}$ satisfies \reff{Normkleinstes} and,
with Proposition~\ref{odd-s1},
\bee
  -1=-\norm{\ti a}_\infty = -\min_{a\in\gG{f}{A}{\conv}}\norm{a}_\infty = 
  \min_{\norm{h}_1\le 1}\gdD{f}{}{A}{h}\,.
\ee
Taking $\ti a=(1,1)\in\gG{f}{A}{\conv}$ we obtain 
\reff{NormkleinstesDual} e.g. for $\ti h=(0,-1)$. However 
$f$ is strictly increasing in the directions $\pm\ti h$ and
$\gdD{f}{}{A}{\ti h}=1$. Hence $\ti h$ is not a descent direction
and \reff{NormkleinstesDual} is not sufficient for their selection. 
Obviously $\ti h=(-1,0)$ is an optimal descent direction on $A$ and 
satisfies  (\ref{NormkleinstesDual}) for every 
$\ti a\in\gG{f}{A}{\conv}$. 

\item Let $X=(\R^2,\norm{\cdot}_\infty)$ and, thus, its dual
$X^*=(\R^2,\norm{\cdot}_1)$. We define $f:\R^2\to\R$ by 
\bee
  f(x,y)=\frac12\left(x+y+|x-y|\right) \,.
\ee 
For $x=(0,0)$ and $A=\sset{x}$ we have 
\bee
  \gG{f}{A}{\conv}=\gG{f}{x}{}=\{ (\la,1-\la) \mid \la\in\left[0,1\right] \}
  \,.
\ee
Again any $\ti a \in\gG{f}{A}{\conv}$ satisfies \reff{Normkleinstes} and,
with Proposition \ref{odd-s1},
\bee
  -1=-\norm{\ti a}_1 = -\min_{a\in\gG{f}{A}{\conv}}\norm{a}_1 = 
  \min_{\norm{h}_\infty\le 1}\gdD{f}{}{A}{h}\,.
\ee
With $\ti a=(1,0)\in\gG{f}{A}{\conv}$ and $\ti h=(-1,1)$ we have 
(\ref{NormkleinstesDual}), but in both directions $\pm\ti h$ function $f$
is strictly increasing and $\gdD{f}{}{A}{\ti h}=1$. Hence $\ti h$ is not a
descent direction and also here \reff{NormkleinstesDual} is not sufficient for
their selection. We readily verify that $\ti h=-(1,1)$ is an optimal descent
direction on $A$ and 
satisfies  (\ref{NormkleinstesDual}) for every 
$\ti a\in\gG{f}{A}{\conv}$. 
\el
\end{exam}
As a consequence of Theorem~{\rm \ref{Satz1Umg}} we obtain that 
descent directions are stable.
\begin{cor}[stability of descent directions]
Let the assumptions of Theorem~{\rm \ref{Satz1Umg}} with 
$0\notin \gG{f}{A}{\conv}$ be satisfied, 
let $\ti a$, $\ti h$ be as there, and
let $L$ be the Lipschitz constant
of $f$ on a neighborhood of $A$. Then every $h\in X$ with
$\|h-\ti h\|<\frac{\norm{\ti a}{}}L$ is a descent direction on $A$. 
\end{cor}
\begin{pf}{}
Let $h\in X$ be as in the statement. By \reff{supportgG} there is 
$a\in \gG{f}{A}{\conv}$ such that 
\ba
  \gdD{f}{}{A}{h} 
&=&
  \sk{a}{h}{} \:=\: \dl a,h-\ti h\dr+\dl a,\ti h\dr
  \stackrel{\tx{Prop. \ref{gg-s1}}}{\leq} 
  L\|h-\ti h\|+\gdD{f}{}{A}{\ti h} \nn\\
&<&
  \|\ti a\| + \gdD{f}{}{A}{\ti h} \stackrel{(\ref{Glamxm})}{=}0\,. \nn
\ea
Hence $h$ is a descent direction. 
\end{pf}
The stability of descent directions allows to work with approximations of an
optimal descent direction. 
\begin{cor}[approximation of an optimal descent
direction]\label{Kiwiel-GS-Verfahren-Corollary}\label{odd-s6}
Let $X$ be uniformly convex (or finite dimensional and
strictly convex) and let $X^*$ be strictly convex.
Moreover let $A\subs X$ be nonempty, let 
$f:X\to\R$ be Lipschitz continuous on a neighborhood of
$A$ with $0\notin \gG{f}{A}{\conv}$, and let $\ti a\in\gG{f}{A}{\conv}$ be as
in Theorem {\rm \ref{Satz1Umg}}. 
Then for any $\de\in]0,1[$ there is some $\tau>0$ such that 
for every $a'\in\gG{f}{A}{\conv}$ with 
\begin{equation*}
 \norm{a'}\leq \min_{a\in \gG{f}{A}{\conv}} \norm{a} + \ta \quad
 \big( = \norm{\ti a}+\ta\: \big)
\end{equation*}
the unique $h'\in X$ satisfying 
\begin{equation}\label{odd-s6-2}
  \sk{a'}{h'}{}=-\norm{a'}\qmq{with} \norm{h'}=1
\end{equation}
is a descent direction on $A$ with 
\begin{equation*}
\Big( \max_{a\in \gG{f}{A}{\conv}} \sk{a}{h'}{}\, = \Big) \quad
\gdD{f}{}{A}{h'} < -\de\norm{\ti a}\,.
\end{equation*}
\end{cor}
Recall that uniformly convex Banach spaces are reflexive
(cf. \cite[Theorem II.2.9]{Cioranescu}) and, thus, the results of Theorem
\ref{Satz1Umg} are available in the corollary.

\begin{pf}{ }
The usual dual mapping $j:X^*\setminus\sset{0}\to \set{x\in X}{\norm{x}=1}$
is given by  
\bee
  \dl a,j(a)\dr =\norm{a}\,.
\ee
Hence \reff{odd-s6-2} just means $h'=-j(a')$ and \reff{NormkleinstesDual}
gives $\ti h=-j(\ti a)$ (notice that 
$0\notin \gG{f}{A}{\conv}$). If the assertion would be false, then
there are $\de>0$ and $a'_k\in \gG{f}{A}{\conv}$ such that 
\bn{odd-s6-4}
  \norm{a'_k}\le\norm{\ti a}+\frac{1}{k} \qmq{and} 
  \gdD{f}{}{A}{-j(a'_k)} \ge -\de\norm{\ti a} \qmz{for all} k\in\N\,.
\ee
By \reff{supportgG} there are $a_k\in\gG{f}{A}{\conv}$ with 
$\gdD{f}{}{A}{-j(a'_k)} = \dl a_k,-j(a'_k)\dr$.
Since $X$ is reflexive and $\gG{f}{A}{\conv}$ weak$^*$-compact, 
we have up to a subsequence that
\bee
  a'_k\wsto : a'\in \gG{f}{A}{\conv} \qmq{and} 
  a_k\wsto : a\in \gG{f}{A}{\conv}\,.
\ee
With \reff{Normkleinstes} we obtain
\bee
  \norm{\ti a} \le \norm{a'} \le \liminf_{k\to\infty} \norm{a'_k} 
  \le \limsup_{k\to\infty} \norm{a'_k}
  \stackrel{(\ref{odd-s6-4})}{\le} \norm{\ti a} \,.
\ee
Since $\ti a$ is uniquely determined by \reff{Normkleinstes}, we get 
$a'=\ti a$ and $\norm{a'_k}\to\norm{\ti a}$. 
Uniform convexity (or finite dimension) of $X$ implies $a'_k\to\ti a$.
Reflexivity of $X$ and strict convexity of $X$ and $X^*$ imply continuity of
$j$ (cf. \cite[Prop. II.5.5]{Cioranescu}) and, thus,
\bee
  -\de\norm{\ti a} \stackrel{(\ref{odd-s6-4})}{\le}
  \liminf_{k\to\infty}\gdD{f}{}{A}{-j(a'_k)} =
  \lim_{k\to\infty}\dl a_k,-j(a'_k)\dr = \dl a,\ti h\dr
  \osm{\tx{\reff{supportgG}}}{\le} \gdD{f}{}{A}{\ti h} 
  \osm{\tx{\reff{Glamxm}}}{=} -\norm{\ti a}\,.     
\ee 
But this is a contradiction and the assertion follows. 
\end{pf}
Let us finally demonstrate with a simple but typical example
how the introduced optimal descent direction 
can improve numerical descent methods.  
\begin{exam}
For $X=\R^2$ equipped with the Euclidean norm we consider 
\bee
  f(x_1,x_2):=|x_1|+\al|x_2| \qmq{with}  0<\al < < 1
\ee 
Here steepest descent methods starting from $(x_1,x_2)$ with 
$x_2>>|x_1|$ easily approach (but usually do not reach) the axis $\{x_1=0\}$ after a few steps. 
Then they highly oscillate around that axis, since the gradients
switch between $(\pm1,\al)$. But with a nonsmooth strategy we would choose 
a suitable ball $A=B_\ep(x)$ at an iteration point $x$ near $\{x_1=0\}$. 
If $0\in B_\ep(x)$, then $0\in \gG{f}{B_\ep(x)}{\conv}$ and we either stop the
algorithm or we decrease ``step size'' $\ep$. If otherwise $0\not\in
B_\ep(x)$, then 
\bee
  \gG{f}{B_\ep(x)}{\conv}=\set{(\la,\al) }{  \la\in[-1,1] }\,.  
\ee
Obviously $\ti a=(0,\al)$ has the smallest norm in 
$\gG{f}{B_\ep(x)}{\conv}$ and the corresponding optimal descent direction on
$B_\ep(x)$ according to Theorem~{\rm \ref{Satz1Umg}} is $\ti h=(0,-1)$.
Now a descent step or a line-search in direction $\ti h$ goes quite
directly to the minimizer $(0,0)$. 
\end{exam}


\end{document}